\documentclass[a4paper,11pt]{article}
\usepackage[dvips]{epsfig}
\usepackage{amsmath,amssymb,euscript,graphicx,amsfonts}
\usepackage{enumerate,color}
\usepackage{graphicx}
\usepackage{hyperref}
\usepackage{tikz}

\newtheorem{theorem}{Theorem}[section]
\newtheorem{proposition}[theorem]{Proposition}
\newtheorem{lemma}[theorem]{Lemma}

\newtheorem{remark}[theorem]{Remark}

\newcommand{\Qed}{\rule{2.5mm}{3mm}}

\renewcommand{\P}{\mathcal{P}}
\newcommand{\la}{\langle}
\newcommand{\ra}{\rangle}
\newcommand{\ZZ}{\mathbb{Z}}

\newcounter{case}

\renewcommand{\thecase}{\arabic{case}}

\newcounter{subcase}

\numberwithin{subcase}{case}

\def\PGL{\hbox{\rm PGL}}
\def\GL{\hbox{\rm GL}}\def\SL{\hbox{\rm SL}}
\def\SO{\hbox{\rm SO}}\def\PSO{\hbox{\rm PSO}}
\def\PSL{\hbox{\rm PSL}}
\def\demo{{\bf Proof}\hskip10pt}
\def\di{\bigm|} \def\lg{\langle} \def\rg{\rangle}
\def\qqed{\hfill $\Box$}

\def\B{{\bf B}}
\def\a{\alpha} \def\b{\beta} \def\g{\gamma} \def\d{\delta} \def\e{\varepsilon}
   
\def\th{\theta}  
\def\D{\Delta}  
\def\si{\Sigma} \def\O{\Omega} \def\G{\Gamma}
\def\f{\noindent}

\def\ol1{\overline 1}  
  
 \def\olt{\overline t} \def\ols{\overline s}  \def\olu{\overline u}
\def\oll{\overline l}

 \def\ox{\overline X}

\def\Aut{\hbox{\rm Aut\,}}\def\soc{\hbox{\rm soc}}
\def\mod{\hbox{\rm mod }}

\def\DD{\mathbb{D}}\def\FF{\mathbb{F}}\def\deg{\hbox{\rm deg}}

\usepackage{minitoc}

\begin{document}

\begin{center}
{\bf\large
Hamilton Cycles In Primitive Graphs of Order $2rs$}
\end{center}

\begin{center}
Shaofei Du{\small \footnotemark}, Yao Tian, Hao Yu\\
\medskip
{\it {\small
Capital Normal University,\\ School of Mathematical Sciences,\\
Bejing 100048, People's Republic of China
}}
\end{center}

\footnotetext{Corresponding author: dushf@mail.cnu.edu.cn.}

\vspace*{-15pt}

\begin{abstract}
After long term efforts, it was recently proved in \cite{DKM2} that  except for the Peterson graph, every connected vertex-transitive graph of order $rs$ has a Hamilton cycle, where $r$ and $s$ are primes.
A natural  topic is to solve the  hamiltonian problem  for  connected vertex-transitive graphs of $2rs$. This topic is  quite trivial, as  the problem is still  unsolved  even for that of $r=3$.
 In this paper, it is shown that except for the Coxeter graph,  every connected vertex-transitive
graph of order $2rs$  contains a Hamilton cycle,  provided  the automorphism group acts  primitively on  vertices.
\end{abstract}

\vspace*{-15pt}

\begin{quotation}
{\em Math. Subj. Class: 05C25, 05C45.}
\end{quotation}

\section{Introduction}
\label{sec:intro}
\indent

Throughout this paper graphs are finite, simple and undirected,
and groups are finite.
Given a graph $X$, by  $V(X)$,  $E(X)$ and $\Aut(X)$ we denote   the
vertex set,  the edge set and the automorphism group of $X$, respectively.
A graph $X$ is {\em vertex-} or {\em arc-transitive}
if $\Aut X$ acts transitively on vertices or arcs, respectively.

Given a transitive group $G$ on  $\O$, a subset $B$ of $\O$ is called a  {\it block} of $G$ if, for any $g\in G$, we have   either $B=B^g$ or $B\cap B^g=\emptyset$. Clearly,
 $G$  has blocks $\O$ and $\{\a\}$ for any $\a\in \O$, which  are said to be {\it trivial}. Then $G$  is  said to be  {\em primitive} if it has no nontrivial blocks.
 Moreover,  a vertex-transitive graph $X$ is said to be {\em primitive}  if $\Aut(X)$  is primitive on vertices.

A simple path (resp. cycle)  containing all vertices of a graph is called a {\it Hamilton path} (resp. {\it cycle}) of this  graph.   A graph containing a Hamilton cycle will be sometimes referred as
a {\em hamiltonian graph.}

In 1970,  Lov\'{a}sz asked in \cite{LL} that
\vskip 3mm
{\it
Does	every finite connected vertex-transitive  graph have a Hamilton	path? }
\vskip 3mm
\f Up to now,   this question remains unresolved and no connected vertex-transitive graph without a Hamilton	path is known to exist.  Moreover, only four (families) of
connected vertex-transitive  graphs	on at least three	vertices not having a
Hamilton	cycle	 are	
known, which are Peterson graph, Coxeter graph and triangle-replaced graphs from them.  Since all of these graphs  are not  Cayley graph, we  may ask if  every connected  Cayley graph has a  Hamilton cycle.

It has been shown that  connected vertex-transitive graphs of orders	
$kp$, $k \leq 6$, $10p$ ($p\ge 11$), $p^j$ ($j \leq 5$) and $2p^2$, where	$p$ is	a prime	
 contain a Hamilton path, see \cite{A79, C98,KM08,KMZ12,KS09,DM85,DM87,MP82,MP83,Z15}.
Furthermore, for all of these families, except for the graphs of orders $6p$ and $10p$ and that four exceptions,
 they contain a Hamilton cycle.
   With the  exception of the Peterson graph, Hamilton cycles are also known to exist in connected
vertex-transitive graphs whose automorphism groups contain a transitive subgroup with a cyclic commutator subgroup of prime-power order (see \cite{DGMW98} and also \cite{D83,DM83,DW85}).

So far we know that Cayley graphs of the following groups contain   a Hamilton cycle:    nilpotent 	
groups of odd order,  with cyclic commutator subgroups
 (see \cite{DGMW98,GWM11,GWM14}); dihedral groups of order divisible by $4$  (see \cite{ACD10}); and
 arbitrary $p$-groups (see \cite{W86}).
   A Hamilton path and in some cases even a Hamilton cycle was proved to exist in
cubic Cayley graphs arising from $(2,\,s,\,3)$-generated
groups (see \cite{GKMM12,GKM09,GM07}).

Recently,   Kutnar, Marusic and the  first author  proved that vertex transitive graphs of  order $rs$  have a Hamilton cycle, except for the Peterson graph (see \cite{DKM1,DKM2}). This work takes many years,  because of  a  difficult  case,  which is a primitive graph with automorphism group  $\PSL(2,\,p)$ and a point-stabilizer  $\DD_{p-1}$.
 A natural question is to consider hamiltonian problem for vertex-transitive graphs of order $2rs$.  As mentioned above, some special cases  have been solved such as that of graphs of order $4p$, $6p$, $10p$ and $2p^2$, where $p$ is a prime (Hamilton path or cycle).
 To solve the general case, a necessary step is to deal with  all primitive graphs of such order. The main result of this paper is the following theorem.
\begin{theorem}
\label{main}
  Except for Coxeter graph, every connected vertex-transitive graph of order $2rs$
contains a Hamilton cycle provided  the automorphism group acts  primitively on its vertices, where $r$ and $s$ are primes.
\end{theorem}

After this introductory section, some notations, basic definitions and useful facts will be given in Section 2  and  Theorem ~\ref{main} will be proved in Section 3.
\section{Preliminaries}
\label{sec:pre}
By $\lfloor a\rfloor$ and $\lceil a\rceil$, we denote the smallest integer no less than $a$ and  largest integer no more than $a$, respectively.
For a prime power $q$, a finite field of order $q$
will be denoted by $\mathbb{F}_q$.
Set  $\mathbb{F}_q^*=\mathbb{F}_q\setminus\{0\}$,
$S=\{t^2\di t\in \mathbb{F}_q\}$,  $S^*=S\cap \mathbb{F}_q^* $ and $N=\mathbb{F}_q^*\setminus S^*$. Then the elements in $S$ and $N$ are called to be {\em squares} and {\em non-squares}, respectively.
By $\ZZ_n$ and $\DD_n$ we  denote a cycle group and dihedral group of order $n$, respectively. For a group $G$ and $L\subset G$, by $C_{G}(L)$ and $N_{G}(L)$ we denote the centralizer and normalizer of $L$ in $G$, respectively.
A semi-product of $K$ and $H$ is denoted by $K\rtimes H$, where $K$ is normal. Let $G$ be a group with a normal subgroup $N$, we denote the image of $g\in G$ under the natural homomorphism of $G$ to $G/N$ by $\overline{g}$. For a group $G$ and its subgroup $H$, $[G:H]$ denotes the set of right cosets of $H$ in $G$; $HgH$ denotes the orbit containing $Hg$ under the action of $H$.

Let $G$ act on some set $\Omega$. For some $\a\in \Omega$ and $g\in G$,  set $\a^G=\{\a^g\di\ g\in G\}.$
 For  $\a\in \Omega$, set $H=G_{\a}$. Then the action of $G$ on $\Omega$ is equivalent to its right multiplication action on right cosets $[G:H]$ relative to $H$.
 For  a subset $\D$ of $\O$, by $G_{(\D)}$ and $G_{\{\D\}}$, we denote the pointwise  and setwise stabilizer of $\D$ in $G$, respectively.

In a graph $X$, let $a\in V(X)$ and $B\subset V(X)$, by $d(a,\,B)$ we denote the number of neighbors of $a$ in $B$.
Given $A$, $B\subset V(X)$, if $d(a,\,B)=d(a',\,B)$ for any $a$, $a'\in A$, then we denote $d(a,\,B)$ by $d(A,\,B)$.

\vskip 3mm
In what follows we recall some definitions related to  orbital graphs and
semiregular automorphisms.
\vskip 3mm
Let $G$ be a transitive permutation group on $\O$. Then $G$ induces
 a natural action on $\O\times \O$.
We call the orbits of $G$ on $\O \times \O$ the {\em orbitals} of $G$,
and in particular  the
{\em trivial} orbital is referred to
$\{(\a,\,\a)\mid\a\in \O \}$. The {\em orbital digraph}
$X(G, \,\G )$ relative to an orbital $\G $ is defined to be the directed
graph with vertex set $\O $ and edge set $\G $. Each orbital $\G $ has
an associated {\em paired orbital} $\G '$ defined by
$\G '=\{ (\b ,\, \a )\di (\a ,\, \b )\in \G \},$ and of course,
$\G $ is said to be {\em self-paired} if $\G =\G '$ in which case
 $X(G,\, \G )$ can be viewed as an undirected graph ({\em orbital graph}).
The $G$-arc-transitive graphs
with vertex-set $\O $ are precisely the orbital graphs $X(G,\, \G )$
for the nontrivial self-paired orbitals $\G $.
In addition, take a point $\a \in \O$, the orbits of the stabilizer
$G_{\a}$ on $\O$
are called {\em suborbits} of $G$ relative to $\a $. There is a one-to-one
correspondence between the suborbits and the orbitals of $G$.
Each orbital $\G_i$ corresponds to a suborbit
$\Delta_i=\{ \,\b\in \O \mid (\a,\,\b)\in \G_i \,\}.$
 Conversely, each suborbit $\Delta_i$ corresponds to an orbital
$\G_i=\{ \, (\a,\,\b)^g\mid g\in G,\ \b \in \Delta_i \,\}$.
 A suborbit of $G$ is said to be
{\em self-paired} if the corresponding orbital is self-paired.
Thus we often use $X(G ,\, \Delta _i )$ and
$X(G,\, \Delta _i\cup \Delta _i')$
to denote graphs $X(G,\, \G )$ and $X(G,\,\, \G \cup \G ')$ respectively.

Let $m\geq 1$ and $n\geq 2$ be integers. An automorphism $\rho$
of a graph $X$ is called $(m,\,n)$-{\em semiregular}
(in short, {\em semiregular})
if as a permutation on $V(X)$ it has a cycle decomposition consisting
of $m$ cycles of length $n$.
If $m=1$ then $X$ is called a {\em circulant}; it is
in fact a Cayley graph of a cyclic group of order $n$.
Let $\P$ be the set of orbits of $\rho$, that is,
the orbits of the cyclic subgroup $\langle \rho \rangle$ generated by $\rho$.
We let the {\em quotient graph corresponding to
$\P$} be the graph $X_\P$ whose vertex set
equals $\P$ with $A,\, B \in \P$ adjacent if there
exist vertices $a \in A$ and $b \in B$, such that $a \sim b$ in $X$.
\vskip 3mm

The following four results will be used later.

\begin{proposition}
\label{pro:preseki}
{\rm \cite[p.167]{MS3}}
Let $F$ be a finite field of odd prime order $p$. Then
$$
|(S^*+1) \cap (-S^*)|=\left\{
\begin{array}{ll}
(p-5)/4 & p\equiv {1\pmod 4},\\
(p+1)/4 & p\equiv {3\pmod 4}.
\end{array}\right.
$$
This implies that if $p\equiv {1\pmod 4}$ then
$$|S^*\cap (S^*+1)|=(p-5)/4,\quad
|N\cap (N+1)| = (p-1)/4,\quad
|S^*\cap (N\pm 1)| = (p-1)/4.$$
\end{proposition}

\begin{proposition}
\label{jack}
{\rm \cite[Theorem~6]{BJ78} (Jackson's Theorem)}
Every $2$-connected regular graph of order $n$ and valency at least $n/3$
contains a Hamilton cycle.
\end{proposition}

\begin{proposition}
\label{Cayley}
{\rm \cite[Corollary~3]{ATD}}
If $X$ is a connected Cayley graph of an abelian group of order at least 3, then every edge of $X$ lies in a hamiltonian cycle.
\end{proposition}

\begin{lemma}
\label{cyclelift}
{\rm \cite[Lemma~5]{MP82} }
Let $X$ be a graph admitting
an $(m,\,p)$-semiregular automorphism $\rho$, where $p$ is a prime.
Let $C$ be a cycle of length $m$ in the quotient graph $X_\P$,
where $\P$ is the set of orbits of $\rho$.
Then, the lift of $C$ either contains a cycle of length
$mp$ or it consists of $p$ disjoint $m$-cycles.
In the latter case we have $d(S,\,S') = 1$ for every edge $SS'$ of $C$.
\end{lemma}

\section{Proof~ of ~Theorem~\ref{main}}
 To prove  Theorem~\ref{main},  let $X$ be a connected vertex-transitive graph of order $2rs$, where $r$ and $s$ are primes. Set $G=\Aut(X)$. It has been proved that $X$ contains a Hamilton cycle if $2rs=2p^2$ or $4p$ for a prime $p$, provided $X$ is not the Coxeter graph which is of order $28$.
   Therefore, in what follows we assume that $r\lneqq s$. If $G$ acts 2-transitively on $V(X)$, then $X$ is a complete graph, which contains a $H$-cycle.
 Now   we need to pick up all the primitive groups of degree $\deg(G)=2rs$ of rank at least $3$ from  \cite{GP} (or \cite{LA03}), where $r$ and $s$ are distinct odd primes.
 Let $H$ be a point stabilizer  in $\soc(G)$.
Checking  \cite{GP}, all the possible groups are listed in Tables 1 and 2.

Table 1 gives the these groups with the socle $\PSL(2,\,q)$.   The first two cases
$H=\DD_{q-1}$ and  $H=\DD_{q+1}$ will be dealt with in Subsections 3.1 and 3.2, respectively.
Moreover,   Magma  shows that every vertex-transitive graph is hamiltonian, arising from other  three groups in Table 1.

Table 2 gives these groups whose socle  is a classical simple group which is not $\PSL(2,\,q)$,  an alternating group or a  sporadic simple group.
These groups  will be dealt with in Subsection 3.3.

\begin{table}
  \centering
   \caption {Primitive groups  of degree $2rs$, where  the socle $\PSL(2,q)$}
 \vskip 2mm \begin{tabular}{ccccc}
   \hline

  No. &$\soc(G)$ & $2rs$           & $ Action$            &$Comment$   \\
\hline

   1& $\PSL(2,\,q)$ & $q(q+1)/2$      &   $G_{\alpha}\cap \soc(G) =\DD_{2(q-1)/d}  $               &$d=(2,\,q-1)$, \\ &&&&$G=\PGL(2,\,11)$ \, {\rm for}\, $q=11$         \\
   2 &$\PSL(2,\,q)$ & $q(q-1)/2$      &   $G_{\alpha}\cap \soc(G) =\DD_{2(q+1)/d}  $               &$d=(2,\,q-1)$          \\
   3 &$\PSL(2,47)$ & $2\times 47\times23$      &   $S_4  $               &~          \\
   4 &$\PSL(2,17)$ & $2\times 17\times3$      &   $S_4  $               &~          \\
   5 &$\PSL(2,41)$ & $2\times 41\times7$      &   $A_5  $               &~          \\

   \hline

 \end{tabular}
 \end{table}

\begin{table}
  \centering
   \caption {Primitive groups $G$ of degree $2rs$, where $\soc(G)\ne \PSL(2,q)$}
 \vskip 2mm \begin{tabular}{ccccc}
   \hline

  No. &$\soc(G)$ & $2rs$           & $ Action$            &$Comment$   \\
\hline
1&$\PSL(4,\,q)$&$\frac{q^3-1}{q-1}(q^2+1)$&2-spaces&$q=3$; or $q=5$; or \\
  ~&~&~&~&$q\equiv11,\,29({\rm mod}~30)$ and \\
   ~&~&~&~&$q$ prime and $q\geq 59$\\
2&$\PSL(5,\,q)$&  $\frac{q^5-1}{q-1}(q^2+1) $ & 2-spaces &$q\equiv-1({\rm mod}~10)$, q prime and $q\geq 29$\\
 3 &${\rm P\Omega}^{-}(2m,\,q)$ & $\frac{(q^{m}+1)(q^{m-1}-1)}{q-1}$      &    on ~t.s. $1$-spaces      &$m$~even          \\
 4&${\rm P\Omega}^{+}(2m,\,q)$ & $\frac{(q^{m}-1)(q^{m-1}+1)}{q-1}$      &    on ~t.s. $1$-spaces         &$m$~odd          \\
 5&$\PSL(3,\,5)$  &$2\times 31\times 3$&   on $(1,\,2)$-dim. flags  &$G=\PSL(3,\,5).2$\\
 6    &$A_{c}$ & $\frac {c(c+1)}2$  &  on 2-sets      & $c\ge 5$                    \\
7 &$M_{11}$ & 66     &  $S_{5}$  &                       \\
 8  &$M_{12}$ & 66     &  $M_{10}:2 $      &                   \\
9 &$M_{23}$ & 506     &  $A_{8} $    &                     \\
10 & $J_{1}$ & 266     &  $\PSL(2,\,11)  $                   &    \\
   \hline
\end{tabular}
 \end{table}

\subsection{ $\soc(G)=\PSL(2,\,q)$  and  $H=\DD_{q-1}$}
 Let $G=\PSL(2,\,q)$ and $H=\DD_{q-1}$.   Consider  the action of $G$ on the set  $[G:H]$ of   cosets  of $H$ in $G$, see  row 1  of Table 1.
Then the degree $n=q(q+1)/2=2rs$,  thus  $q\equiv 3({\rm mod} ~4)$ and in particular $-1\in N$, the set of non-squares. Set $\FF_q^*=\lg \th \rg.$

For any $g\in \SL(2,\,q)$, set $\overline{g}=gZ(\SL(2,\,q))$.
In $\SL(2,\,q)$, set
$$u=\left(
             \begin{array}{cc}
               1 & 1 \\
               0 & 1 \\
             \end{array}
           \right),\,
           u'=\left(
             \begin{array}{cc}
               1 & 0 \\
               1 & 1 \\
             \end{array}
           \right),\,
          l=\left(
             \begin{array}{cc}
               \th & 0 \\
               0 & \th^{-1} \\
             \end{array}
           \right),\,
           t=\left(
             \begin{array}{cc}
               0 & 1 \\
               -1 & 0 \\
             \end{array}
           \right).$$
           Since $\PSL(2,\,q)$ has only one conjugacy class of subgroups isomorphic to $\DD_{q-1}$, we may set $H=\lg \oll,\, \olt \rg $.
Let $V$ be the row vector space so that  the action of $g\in \GL(2,\,q)$ on a vector $(x,\,y)$ is just defined as $(x,\,y)\cdot g$.  Set $\frac yx=\lg (x,\,y)\rg $. Then   all the  projective points are  $\{\infty,\,0,\,1,\,2,\,\cdots,\,q-1\}$.
The action of $G$ on $[G:H]$ is equivalent to its  action  on the set of unordered pairs of distinct projective points, where $H=G_{ \{0,\, \infty\}}$. Thus we have
$$\begin{array}{ll} &\olu ':   \{\infty,\, 0\}\to \{1,\, 0\},\quad
\overline{l^i}:
  \{j,\, j+1\}\to \{j\th^{-2i},\, (j+1)\th^{-2i}\}, \\ &\overline{l^it}: \{j,\, j+1\}\to \{-j^{-1}\th^{2i},\, -(j+1)^{-1}\th^{2i}\}.\end{array}$$
Then the all $\lg \olu \rg$-orbits  are
$$B_\infty=\{\{\infty,\,i\}|i\in \FF_q\},\quad B_j=\{\{i,\,i+j\}|i\in \FF_q\},\,\, j\in\{ 1,\,2,\,3,\,\cdots,\,\frac{q-1}{2}\}.$$
  Set ${\bf B}=\{ B_j\di j\in \{ 1,\,2,\,3,\,\cdots,\,\frac{q-1}{2}\}\}$.
Considering  the action of $N_G(\lg \olu\rg )=\lg \olu\rg \rtimes \lg \oll\rg$  on the vertices, we know that
  $N_G(\lg \olu\rg)$ fixes the  block $B_\infty$  setwise and acts transitively on other vertices. In particular, $\lg \oll\rg$ fixes $B_\infty$ and acts  regularly on  $\frac{q-1}2$ remaining blocks $B_j$ in ${\bf B}$.

The suborbits of $G$ have been determined in  \cite{LWWX94} and  an alternating  description  is given below.

\vskip 3mm
\begin{lemma}\label{suborbit}  Suppose $q\equiv 3(\mod 4)$. Then every nontrivial  suborbit of $G$ relative to $H$  can be written as  $\{j,\,j+1\}^H$, where $j\in \FF_q$, with length  $\frac{q-1}2$ and $q-1$ if and only if $j^{2}+j \in N$ and $j^2+j\in S$, respectively. Moreover,
  $\{j,\,j+1\}^H$ is self-paired if and only if either $j+1\in N$ or $j\in S$, and if it is non self-paired, then its paired suborbit is    $\{-j,\,-j-1\}^H$.
 \end{lemma}
 \demo For $i\in \FF^*_q$, direct computations show that $\{\infty, i\}$ belongs to $\{0,\,1\}^H$ or $\{0,\,-1\}^H$ depending on whether $i\in S^*$ or $i\in N$, respectively. Since $\lg \oll\rg\leq H$ acts regularly on ${\bf B}$, any other suborbits can also be written as $\{j,\,j+1\}^H$. The length of $\{j,\,j+1\}^H$ is $\frac{q-1}2$ and $q-1$ if and only if the order of the stabilizer for $\{j,\,j+1\}$ in $H$ is $2$ and $1$, respectively. But the former holds   if and only if there exists some $k\in \ZZ_q$ such that $\overline{l^kt}$ fixes $\{j,\,j+1\}$, i.e., $j+1=-j^{-1}\th^{2k}$. Therefore we deduce that the length of the suborbit is  $\frac{q-1}2$ or $q-1$ depending on  $j^{2}+j \in N$ or $j^2+j\in S$, respectively.

Let $\D=\{j,\, j+1\}^H$. If $j+1=0$, then  $\D^*=\{0,\,1\}^H$. If $j+1\neq0$, then  $\D^*=\{\frac{-j}{j+1},\,-1\}^H$. Now, $\D$ is self-paired if and only if there exists some element of $H$
 mapping  $\{j,\, j+1\}$ to $\{\frac{-j}{j+1},\,-1\}$. From
 $$\{\overline{l^k}(j),\, \overline{l^k}(j+1)\}=\{j\th^{-2k}, \,(j+1)\th^{-2k}\}=\{\frac{-j}{j+1},\,-1\}\quad {\rm and}$$
  $$\{\overline{l^kt}(j),\, \overline{l^kt}(j+1)\}=\{-j^{-1}\th^{2k}, \, -(j+1)^{-1}\th^{2k}\}=\{\frac{-j}{j+1},\,-1\},$$
 we  know that such element  of $H$ exists if and only if  $j+1\in N$ or $j\in S$, as desired.

 Suppose that  $\D$ is not self-paired and $j+1\neq0$. Then $j+1=\th^{-2k}\in S$ and  $\overline{l^k}$ maps $\{\frac{-j}{j+1},\,-1\}$ to $\{-j,\,-j-1\}$, that is $\D^*=\{\frac{-j}{j+1},\,-1\}^H=\{-j,\,-j-1\}^H.$\qqed

\begin{remark} By Lemma~\ref{suborbit}, it is easy to determine the number of nontrivial suborbits of length $\frac{q-1}2$ or $q-1$,  and  the number of nontrivial paired suborbits. But we do not need these numbers in here.
\end{remark}

Before going to prove the main result,   we first  give a technical lemma on number theory.

\begin{lemma}\label{feng}
Suppose that $q$ is an odd prime power. If $a,\,b\in \FF_q^*$ and $a\neq b$. Then
$$\begin{array}{ll} &|(S^*+a)\cap (S^*+b)\cap N|\leq \lceil\frac{1}{8}(q+11+2\sqrt{q})\rceil ,\\
 &|(S^*+a)\cap (N+b)\cap N|\leq \lceil\frac{1}8(q+11+2\sqrt{q})\rceil ,\\
&|(S^*+a)\cap (N+b)\cap S^*|\geq\lfloor\frac{1}{8}(q-11-2\sqrt{q})\rfloor, \\
&|(N+a)\cap (N+b)\cap S^*|\geq\lfloor\frac{1}{8}(q-11-2\sqrt{q})\rfloor. \end{array} $$

\end{lemma}
\demo
Set $\eta: \FF^*_q\rightarrow\{\pm1\}$ by assigning the elements of $S^*$ to $1$ and that of $N$ to $-1$ and moreover, set $\eta(0)=0$. This $\eta$ is exactly that in \cite[Example 5.10]{LN}. Also we need to quote the following three results from~\cite[Theorems 5.4, 5.48, 5.41]{LN}:
\begin{enumerate}
\item[\rm(i)] $\sum_{x\in \FF_q}\eta(x)=0$;
\item[\rm(ii)] $\sum_{x\in \FF_q}\eta(x^2+Ax+B)=q-1$ for $A^2-4B=0$ or $-1$ for otherwise, where $A,B\in \FF_q$;
\item[\rm(iii)] $|m|\leq 2\sqrt{q}$, where  $m:=\sum_{x\in \FF_q}\eta(x(x-1)(x-t))$ and $t\in \FF_q$.
\end{enumerate}

For four  inequalities of the lemma, we have the same arguments and  here we  just prove the first one.  Set $W=(S^*+a)\cap (S^*+b)\cap N$, that is $$W=\{x\in\FF_q\mid \eta(x-a)=\eta(x-b)=1,\, \eta(x)=-1\}.$$
Now let $a,\,b\in S^*$.
Then by the above three formulas (i)-(iii), we have
$$\begin{array}{lcl}
 |W|&=&\frac{1}{8}\sum_{x\in \FF_q\setminus\{0,a,b\}}(1+\eta(x-a))(1+\eta(x-b))(1-\eta(x))\\
     &=&\frac{1}{8}\sum_{x\in \FF_q\setminus\{0,a,b\}}(1-\eta(x)+\eta(x-a)+\eta(x-b)-\eta(x(x-a))-\eta(x(x-b))\\
     &&+\eta((x-a)(x-b))-\eta((x-a)(x-b)x)\\
     &=&\frac{1}{8}[(q-3)-(-\eta(b)-\eta(a))-(\eta(-a)+\eta(b-a))-(\eta(-b)+\eta(a-b))\\
     &&-(-1-\eta b(b-a))-(-1-\eta a(a-b))+(-1-\eta(ab))+m]\\
     &\leq&\lceil\frac{1}{8}(q+11+2\sqrt{q})\rceil.
\end{array}$$
\qqed

According to Lemma~\ref{suborbit}, we shall deal with the orbital graphs $X=X(G,\,\D)$ or $X=X(G,\,\D\cup \D^*)$, according to that  $\D$ is  self-paired and of length $\frac{q-1}2$,  non self-paired and of length $\frac{q-1}2$,
 self-paired and of length $q-1$,  and non self-paired and of length $q-1$, respectively, in the following four lemmas.

\begin{lemma}\label{self} Suppose that $\D$ is  self-paired and of length $\frac{q-1}2$. Then
$X(G,\,\D)$  is hamiltonian.
\end{lemma}
\demo Let $X=X(G,\,\D)$,  where $\D$ is  self-paired and of length $\frac{q-1}2$. Let $Y$ be the quotient graph induced by $\lg \olu\rg $, with vertices $\B\cup \{B_\infty\}$. Then by
 Lemma~\ref{suborbit}, we may set $\D=\{j,\, j+1\}^H$, where $j(j+1)\in N$, $ j+1\in N$ and $j\in \FF_q.$
Then  the neighborhood of $\{0,\infty\}$ is:
$$X_1(\{0,\, \infty\})=\D=\{ \{ j\th^{-2k} ,\, (j+1)\th ^{-2k}\}\di k\in \FF_q\}.$$
 Since $|\D|=\frac{q-1}2$ and $\lg \oll\rg$  acts regularly on ${\bf B}$,   $d(B_\infty,\, B_i)=1$ for any $i=1,$ $2,$ $3,$ $\cdots,$ $\frac{q-1}{2}$.

 \vskip 3mm
 The lemma will be proved by the following three steps:

 \vskip 3mm
 {\it Step 1: Show  $d(B_m,\,B_i)\le 2$ for any $i,$ $ m=1,$ $2,$ $3,$ $\cdots,$ $\frac{q-1}{2}$.}
\vskip 3mm

 Since $\lg \overline{r}\rg $ is regular on ${\bf B}$ and $\{0,\, 1\}\in B_1$,  we may just consider $d(B_1,\, B_i)=d(\{0,\, 1\},\,B_i)$ for any $i=1,$ $2,$ $3,$ $\cdots,$ $\frac{q-1}{2}$.
  Since $\olu '$ maps $\{\infty, \,0\}$ to $\{0, \,1\}$, we know that
 $$X_1(\{0,\,1\})=\D^{\olu'}=\{ \{ j\th^{-2k},\, (j+1)\th ^{-2k}\}\di k\in \FF_q\}^{\olu'}=\{ \{\frac{j\th^{-2k}}{1+j\th^{-2k}},\,\frac{(j+1)\th^{-2k}}{1+(j+1)\th^{-2k}}\}\di k\in \FF_q \}.$$
So a vertex in $X_1(\{0,\,1\})$   is contained  in $B_i$ if and only if
$$\{\frac{j\th^{-2k}}{1+j\th^{-2k}},\,\frac{(j+1)\th^{-2k}}{1+(j+1)\th^{-2k}}\}=\{t,\,t+i\}~{\rm for~ some}~t,$$
if and only if  one  of the following two systems of equations  has solutions:
\begin{eqnarray}
\label{n4}\frac{j\th^{-2k}}{1+j\th^{-2k}}=t,\quad \frac{(j+1)\th^{-2k}}{1+(j+1)\th^{-2k}}=t+i;
\end{eqnarray} and
\begin{eqnarray}
\label{n5}\frac{j\th^{-2k}}{1+j\th^{-2k}}=t+i,\quad \frac{(j+1)\th^{-2k}}{1+(j+1)\th^{-2k}}=t.
\end{eqnarray}

Solving Eq(\ref{n4}),   we  get $$ij(j+1)u^2+(2ij+i-1)u+i=0,$$
where $u=\th^{-2k}$.   This  equation has solutions for $u$  if and   only if
$$\d_1:=(2ij+i-1)^2-4i^2j(j+1)=i^2-(2+4j)i+1\in S^*.$$
 Since the product of two solutions $u_1$ and $u_2$  is $(j(j+1))^{-1}$, a non-square,   we know that either $u_1\in S^*$  or $u_2\in S^*$ if the above equation has solutions. Therefore,  there exists exactly one solution  for $\th^{-2k}=u$   if and only if $\d_1\in S^*$, noting  that every $\th^{-2k}$ gives a unique $t$, equivalently, a unique vertex in the block $B_i$.

Solving Eq(\ref{n5}),   we  get $$ij(j+1)u^2+(2ij+i+1)u+i=0,$$
where $u=\th^{-2k}$.   This  equation has solutions for $u$  if and   only if
$$\d_2:=(2ij+i+1)^2-4i^2j(j+1)=i^2+(2+4j)i+1\in S^*.$$
Similarly,    there exists exactly   one  solution  for $\th^{-2k}$   if and only if $\d_2\in S^*$.

\vskip 3mm
Summarizing   Eq(\ref{n4}) and  Eq(\ref{n5}), we get $d(\{0, \,1\},\, B_i)\le 2$.

 \vskip 3mm
 {\it Step 2: Show  that for a given $j$,  there exists some $i$ such that $d(B_j,\,B_i)=2$.}
\vskip 3mm

It suffices to show
   $d(\{0,\,1\},\, B_i)=2$ for some $i\ne 0$, equivalently,   to show that  the number of  $B_i$ ($i\ne 1$)   such that   $d(B_1,\,B_{i})=1$ is  less than   $\frac{q-1}2-1-2=\frac{q-7}2$.

Now, $d(B_1,\, B_i)=1$  if and only if
$$\d_1\d_2=(i^2-(2+4j)i+1)(i^2+(2+4j)i+1)=y\in N,$$
that is
\begin{eqnarray}
\label{n6}u^2+(2-(2+4j)^2)u+1-y=0,\end{eqnarray}
where $u=i^{2}$. Note that for a given $u\in S^*$, $i$ and  $-i$ give the same block $B_i$. Thus a solution of $u$ can provide at most one block $B_i$ satisfying our conditions.

 In what follows, we analyse the number of solutions for $u$.

 Eq(\ref{n6})  has some solutions for $u$ if and only if
$$\d:=(2-(2+4j)^2)^2-4(1-y)\in S,$$
that is
$$y\in S+t,\, {\rm where}\,\, t=-4j(j+1)\in S.$$
Now $y\in (S+t)\cap N$.
First suppose that $1-y\in N$. Then $y\in (S+t)\cap N\cap (1+S)$. By Lemma \ref{feng}, we have at most $\lceil\frac{1}{8}(q+11+2\sqrt{q})\rceil+1$ choices for $y$, and then for $u$ as well.
Secondly, suppose that $1-y\in S$. Then $y\in (S+t)\cap N\cap (1+N)$. By Lemma \ref{feng}, we have at most $\lceil\frac{1}{8}(q+11+2\sqrt{q})\rceil+1$ choices for $y$. Since every $y$ may give two solutions for $u$, we have  at most $2\lceil\frac{1}{8}(q+11+2\sqrt{q})\rceil+2$ solutions for $u$.

In summary, we have at most
$$ \lceil\frac{1}{8}( q+11+2\sqrt{q})\rceil+2\lceil\frac{1}{8}( q+11+2\sqrt{q})\rceil+3$$
\f  blocks $B_i$ such that $d(B_0,\, B_i)=1$. Now
$$\lceil\frac{1}{8}( q+11+2\sqrt{q})\rceil+2\lceil\frac{1}{8}( q+11+2\sqrt{q})\rceil+3\le \frac{q-7}2,$$
 provided  $q>169$.  In other words,  if $q>169$ there exists some $i$ such that $d(B_0,\,B_i)=2$. For $7\le q\le 169$, only the primes 19, 43, 67 and 163 satisfy $\frac{q(q+1)}2=2rs$. For these primes, we can get a Hamilton cycle by Magma.

\vskip 3mm
{\it Step 3: Show the existence of a $H$-cycle.}
\vskip 3mm
Let's come back to the proof of the lemma.  Let $Y_1=Y[{\bf B}]$, the subgraph of $Y$ induced by ${\bf B}$. Then $Y_1$ is a Cayley graph on $\ZZ_{\frac {q-1}2}$. Since the valency of $X$ is $\frac{q-1}2$, $d(B_1,\, B_\infty)=1$,
and $d(B_1,\, B_i)\le 2$, it  follows from
$$\frac 12(\frac{q-1}2-1-2)\ge \frac 13\cdot \frac{q-1}2$$
that $Y_1$ has at most two connected components. Since $\frac{q-1}2$ is odd,  $Y_1$  must be connected.
 Now there are double edges between $B_1$ and $B_i$ for some $i$.
By Proposition~\ref{Cayley}, $Y_1$ contains a cycle  passing the edge $B_1 B_i$,   say $\cdots B_j B_1 B_i \cdots $.   In $Y$,  replacing  the edge $B_j B_1$  by the path $B_j B_\infty B_1$,
we get a $H$-cycle, say $C$ for $Y$. By Proposition~\ref{cyclelift}, $C$ can be lifted to a $H$-cycle of $X$.  \qqed

\begin{lemma}
Suppose that $\D$ is non self-paired and of length $\frac{q-1}2$. Then
$X(G,\,\D)$  is hamiltonian.
\end{lemma}
\demo
Let $X=X(G,\,\D\cup\D^*)$,  where $\D$ is non self-paired and of length $\frac{q-1}2$. Let $Y$ be the quotient graph induced by $\B\cup \{B_\infty\}$. Then by
 Lemma~\ref{suborbit}, we may set $\D=\{j,\, j+1\}^H$ and $\D^*=\{-j,\, -j-1\}^H$ where $j(j+1)\in N$, $j+1\in S$,\, $j \in N$ and $j\in \FF_q.$
Then  the neighborhood of $\{0,\,\infty\}$ is:
$$X_1(\{0,\, \infty\})=\D\cup\D^*=\{ \{ j\th^{-2k} ,\, (j+1)\th ^{-2k}\},\,\{ (-j)\th^{-2k},\, (-j-1)\th ^{-2k}\}\di k\in \FF_q\}.$$
 Since $|\D\cup\D^*|=q-1$ and $\lg \oll\rg$  acts regularly on ${\bf B}$, $d(B_\infty,\, B_i)=2$ for any $i= 1,\,2,\,\cdots,\,\frac{q-1}{2}$.
\vskip 3mm
 The lemma will be proved by the following two steps:

\vskip 3mm
{\it Step 1:  $d(B_k,\,B_i)\in \{0,\,2,\,4\}$ for any $i,\,k=1,\,2,\,\cdots,\,\frac{q-1}{2}$.}
    \vskip 3mm
Since $\lg \oll\rg $ is regular on ${\bf B}$ and $\{0,\,1\}\in B_1$,  we may just consider $d(B_1,\,B_i)=d(\{0, 1\},\,B_i)$ for any $i=1,\,2,\,3,\cdots,\,\frac{q-1}{2}$. Since $\olu '$ maps $\{\infty,\,0\}$ to $\{0,\,1\}$,
 we know that
 $$X_1(\{0,\,1\})=\{\D,\,\D^*\}^{\olu '}=\{ \{\frac{j\th^{-2k}}{1+j\th^{-2k}},\,\frac{(j+1)\th^{-2k}}{1+(j+1)\th^{-2k}}\},\, \{\frac{(-j)\th^{-2k}}{1+(-j)\th^{-2k}},\,\frac{(-j-1)\th^{-2k}}{1+(-j-1)\th^{-2k}}\}\di k\in \FF_q \}.$$
A vertex in $X_1(\{0,\,1\})$   is contained  in $B_i$ if and only if  some of the following four systems of equations has solutions:
\begin{eqnarray}
\label{n7}\frac{j\th^{-2k}}{1+j\th^{-2k}}=t,\quad \frac{(j+1)\th^{-2k}}{1+(j+1)\th^{-2k}}=t+i;
\end{eqnarray}
\begin{eqnarray}
\label{n8}\frac{j\th^{-2k}}{1+j\th^{-2k}}=t+i,\quad \frac{(j+1)\th^{-2k}}{1+(j+1)\th^{-2k}}=t;
\end{eqnarray}
\begin{eqnarray}
\label{n9}\frac{(-j)\th^{-2k}}{1+(-j)\th^{-2k}}=t,\quad \frac{(-j-1)\th^{-2k}}{1+(-j-1)\th^{-2k}}=t+i;
\end{eqnarray}
\begin{eqnarray}
\label{n10}\frac{(-j)\th^{-2k}}{1+(-j)\th^{-2k}}=t+i,\quad \frac{(-j-1)\th^{-2k}}{1+(-j-1)\th^{-2k}}=t.
\end{eqnarray}

Solving Eq(\ref{n7}) and Eq(\ref{n9}),   we  get the respective equation $$ij(j+1)u^2\pm (2ij+i-1)u+i=0,$$
where $u=\th^{-2k}$.   For each of these two  equations, it  has solutions for $u$  if and   only if
$$\d_1:=(2ij+i-1)^2-4i^2j(j+1)=i^2-(2+4j)i+1\in S^*.$$
 Since the product of two solutions $u_1$ and $u_2$  is $(j(j+1))^{-1}$, a non-square,   we know that either $u_1\in S^*$  or $u_2\in S^*$ if the above equation has solutions. Therefore,  there exists exactly one solution  for $\th^{-2k}=u$   if and only if $\d_1\in S^*$, noting  that every $\th^{-2k}$ gives a unique $t$, equivalently, a unique vertex in the block $B_i$. Totally, two systems of equations give two vertices in the $B_i$.

Solving Eq(\ref{n8}) and Eq(\ref{n10}),   we  get respective  equation $$ij(j+1)u^2\pm (2ij+i+1)u+i=0,$$
where $u=\th^{-2k}$.   This  equation has solutions for $u$  if and   only if
$$\d_2:=(2ij+i+1)^2-4i^2j(j+1)=i^2+(2+4j)i+1\in S^*.$$
Similarly,    there exists exactly   one  solution  for $\th^{-2k}$   if and only if $\d_2\in S^*$. Totally, two systems of equations give two vertices in the $B_i$.

\vskip 3mm
In summary,  $d(B_1,\,B_i)=2$ if and only if $\d_1\d_2\in N$;  and  $d(B_1,\,B_i)=0$ or $4$ provided $\d_1\d_2\in S$.
\vskip 3mm
{\it Step 2: Show the existence of a $H$-cycle.}
\vskip 3mm

Let $Y_1=Y[{\bf B}]$, the subgraph of $Y$ induced by ${\bf B}$. Then $Y_1$ is a Cayley graph on $\ZZ_{\frac {q-1}2}$. Since the valency of $X$ is $q-1$, $d(B_1,\,B_\infty)=2$,
and $d(B_1,\,B_i)\le 4$, it  follows from
$$\frac 14(q-1-4-2)\ge \frac 13\cdot \frac{q-1}2$$
that $Y_1$ has at most two connected components. Then, using the  same arguments in Step 3 of  Lemma~\ref{self},  one may get a $H$-cycle of $X$. \qqed

\begin{lemma}
Suppose that $\D$ is  self-paired and of length $q-1$. Then
$X(G,\,\D)$  is hamiltonian.
\end{lemma}
\demo
Let $X=X(G,\,\D)$,  where $\D$ is self-paired and of length $q-1$. Let $Y$ be the quotient graph induced by $\B\cup \{B_\infty\}$. Then by
 Lemma~\ref{suborbit}, we may set $\D=\{j,\,j+1\}^H$ where $j(j+1)\in S^*$ and either $j+1\in N$ or $j\in S^*.$
Then  the neighborhood of $\{0,\,\infty\}$ is:
$$X_1(\{0,\, \infty\})=\D=\{ \{ j\th^{-2k},\, (j+1)\th ^{-2k}\},\,\{ (-j)\th^{-2k},\, (-j-1)\th ^{-2k}\}\di k\in \FF_q\}.$$
 Since $|\D|=q-1$ and $\lg \overline{l}\rg$  acts regularly on ${\bf B}$,\,$d(B_\infty, B_i)=2$ for any $i= 1,\,2,\,\cdots,\,\frac{q-1}{2}$.
\vskip 3mm
 The lemma will be proved by the following two steps:

\vskip 3mm
{\it Step 1:  $d(B_m,\,B_i)\le 4$ for any $i,\,m\in \FF_q^*$.}
    \vskip 3mm
Since $\lg \overline{l}\rg $ is regular on ${\bf B}$ and  $\{0,\,1\}\in B_1$,  we may just consider $d(B_1,\, B_i)=d(\{0,\, 1\},\,B_i)$ for any $i\in \FF_q^*$.
Now,
 $$X_1(\{0,\,1\})=\{\D\}^{\olu '}=\{ \{\frac{j\th^{-2k}}{1+j\th^{-2k}},\,\frac{(j+1)\th^{-2k}}{1+(j+1)\th^{-2k}}\}, \, \{\frac{(-j)\th^{-2k}}{1+(-j)\th^{-2k}},\,\frac{(-j-1)\th^{-2k}}{1+(-j-1)\th^{-2k}}\}\di k\in \FF_q \}.$$
A vertex in $X_1(\{0,\,1\})$  is contained  in $B_i$ if and only if one of the following four systems of equations has solutions:
\begin{eqnarray}
\label{n11}\frac{j\th^{-2k}}{1+j\th^{-2k}}=t,\quad \frac{(j+1)\th^{-2k}}{1+(j+1)\th^{-2k}}=t+i;
\end{eqnarray}
\begin{eqnarray}
\label{n12}\frac{j\th^{-2k}}{1+j\th^{-2k}}=t+i,\quad \frac{(j+1)\th^{-2k}}{1+(j+1)\th^{-2k}}=t;
\end{eqnarray}
\begin{eqnarray}
\label{n13}\frac{(-j)\th^{-2k}}{1+(-j)\th^{-2k}}=t,\quad \frac{(-j-1)\th^{-2k}}{1+(-j-1)\th^{-2k}}=t+i;
\end{eqnarray}
\begin{eqnarray}
\label{n14}\frac{(-j)\th^{-2k}}{1+(-j)\th^{-2k}}=t+i,\quad \frac{(-j-1)\th^{-2k}}{1+(-j-1)\th^{-2k}}=t.
\end{eqnarray}

Solving Eq(\ref{n11}) and  Eq(\ref{n13})    we  get the respective equation  $$ij(j+1)u^2\pm (2ij+i-1)u+i=0,$$
 where $u=\th^{-2k}$.   For each of these two equations has solutions for $u$ only if
$$\d_1:=(2ij+i-1)^2-4i^2j(j+1)=i^2-(2+4j)i+1\in S.$$
(1)$\d_1\in S^*$: Since the product of two solutions $u_1$ and $u_2$  is $(j(j+1))^{-1}$, a square,   we know that either $u_1,\,u_2\in S^*$  or $u_1,\, u_2\in N^*$. Therefore,  there exist two solutions  for $\th^{-2k}=u$   only if $\d_1\in S^*$. Noting  that every $\th^{-2k}$ gives a unique $t$, equivalently, one vertex in the block $B_i$.  Thus two systems of equations  give   two vertices in $B_i$.\\
(2)$\d_1=0$: For these two equations, there is just one solution for $u$ and it gives a unique $t$.   Thus two systems of equations  give   one vertice in $B_i$.

Solving Eq(\ref{n12}) and  Eq(\ref{n14}), we  get respective equation  $$ij(j+1)u^2\pm (2ij+i+1)u+i=0,$$
where $u=\th^{-2k}$.   This  equation has solutions for $u$  if and   only if
$$\d_2:=(2ij+i+1)^2-4i^2j(j+1)=i^2+(2+4j)i+1\in S.$$
Similarly, if $\d_2 \in S^*$,  there exist exactly  two  solutions  for $\th^{-2k}$.  Thus two equations  give   two vertices in $B_i$.  If $\d_2 =0$, there exists one  solution  for $\th^{-2k}$. Thus we only get one vertex in $B_i$.

In summary,  $d(B_1,\,B_i)=2$ if and only if $\d_1\d_2\in N$; $d(B_1,\,B_i)=0$ or $4$, provided  $\d_1\d_2\in S^*$; and $d(B_1,\,B_i)=1 ~{\rm or}~3$ if and only if $\d_1\d_2=0$.

\vskip 3mm
{\it Step 2: Show the existence of a $H$-cycle.}
\vskip 3mm

Let $Y_1=Y[{\bf B}]$ be the subgraph of $Y$ induced by ${\bf B}$. Then $Y_1$ is a Cayley graph on $\ZZ_{\frac {q-1}2}$. Since the valency of $X$ is $q-1$, $d(B_1,\, B_\infty)=2$,
and $d(B_1,\, B_i)\le 4$, it  follows from
$$\frac 14(q-1-4-2)\ge \frac 13\cdot \frac{q-1}2.$$
Then we get a $H$-cycle, with the same arguments as in Step 3 of Lemma~\ref{self}.
\qqed

\begin{lemma}
Suppose that $\D$ is non self-paired and of length $q-1$. Then
$X(G,\,\D)$  is hamiltonian.
\end{lemma}
\demo In this case, $\D=\{1,\,0\}^H$ and $\D^*=\{-1,\,0\}^H$. Let $X=X(\D\cup\D^*)$ and $Y$ the quotient graph induced by $\B\cup \{B_\infty\}$.  Then  the neighborhood of $\{0,\,\infty\}$ is:
$$X_1(\{0,\, \infty\})=\D\cup \D^*=\{ \{ 0,\,\th^k\},\,\{ \infty,\,\th^k\}\di k\in \FF_q\}.$$
By observing the vertices of block $B_\infty$, we get $d(B_\infty)=q-1$, and since $\lg \overline{l}\rg$ is regular on ${\bf B}$,  $d(B_\infty,\, B_i)=2$ for any $i=1,\,2,\,\cdots,\,\frac{q-1}{2}$.  Since $\olu '$ maps $\{\infty,\, 0\}$ to $\{0,\, 1\}$, we know that
 $$X_1(\{0,\,1\})=\{\D,\,\D^*\}^{\olu '}=
\{\{0,\,\frac{\th^k}{1+\th^k}\},\,\{1,\,\frac{\th^k}{1+\th^k}\}\di k\in \FF_q\}.$$

A direct computation shows $d(B_1)=2$. Moreover,  $d(B_1,\,B_i)$ is exactly   the number of union of solutions of the following two equations:
 $$\{0,\,\frac{\th^k}{1+\th^k}\}=\{t+i,\,t\}\quad {\rm and} \quad \{1,\,\frac{\th^k}{1+\th^k}\}=\{t+i,\,t\}.$$
Solving them, we get four solutions:
 $$\th^k=\frac{-i}{1+i},\,t=-i; \quad \th^k=\frac{i}{1-i},\,t=0; \quad\th^k=\frac{1-i}i,\,t=1-i;\quad  {\rm and}\,~  \th^k=\frac{-i-1}{i},\,t=1.$$
Therefore, $d( B_1,\,B_i)=4.$

Since $\lg \overline{l} \rg $ is regular on ${\bf B}$,  $d(B_j,\,B_i)=d(B_1,\,B_{i'})$ for some $i'$ and $d(B_i)=d(B_1)$. Then we conclude that $d(B_i,\,B_j)=4$ and $d(B_i)=2$. Thus the graph $Y\setminus \{B_\infty\}$ is a complete graph.  As before,   $X$ is hamiltonian.\qqed

\subsection {$\soc(G)=\PSL(2,\,q)$ and $H=\DD_{q+1}$}
Let $G=\PSL(2,\,q)$ and $H=\DD_{q+1}$.
 Consider  the action of $G$ on the set  $[G:H]$ of   cosets  of $H$ in $G$, see  row 2  of Table 1.
 Then $n=\frac{q(q-1)}2=2rs$. This implies that  $q\equiv 1({\rm mod}~ 4)$ and both  $q$ and $\frac{q-1}{4}$ are primes. So $r=\frac{q-1}{4}$ and $s=q$.
 Set $\FF_q^*=\lg \th \rg $ and  $\sqrt{-1}=\th^{\frac{q-1}4}$. In $\GL(2,\,q)$, we set
$$u=\left(
             \begin{array}{cc}
               1 & 1 \\
               0 & 1 \\
             \end{array}
           \right),\,
           u'=\left(
             \begin{array}{cc}
               1 & 0 \\
               1 & 1 \\
             \end{array}
           \right),\,
          l=\left(
             \begin{array}{cc}
               \th & 0 \\
               0 & \th^{-1} \\
             \end{array}
           \right),\,
           t=\left(
             \begin{array}{cc}
               0 & 1 \\
               -1 & 0 \\
             \end{array}
           \right),$$
          $$t(x,\,y)=\left(
             \begin{array}{cc}
               x & y\th \\
               y & x \\
             \end{array}
           \right),\,
           t'(x,\,y)=\sqrt{-1}\left(
                              \begin{array}{cc}
                                1 & 0 \\
                                0 & -1 \\
                              \end{array}
                            \right)t(x,\,y)
           =\sqrt{-1}\left(
             \begin{array}{cc}
               x & -y\th \\
               y & -x \\
             \end{array}
           \right),\,x\neq 0.$$
  Then up to conjugacy,
  $H$ may be chosen as
$$H=\{\overline{t(x,\,y)}, \overline{t'(x,\,y)}\di x^2-y^2\th=1\}.$$

Consider the action of $N_G(\lg \olu\rg)=\lg \olu\rg \rtimes \lg \oll\rg$ on the set of $\lg \olu \rg$-orbits (blocks) on $[G:H]$. Then $[G:H]$  can be divided into two parts, say $\B$ and $\B'$,  where
 $$ \B=\{B_{1},\,B_{2},\,\cdots,\,B_{\frac{q-1}4}\}, \quad \B'=\{B'_{1},\,B_{2}',\,\cdots,\,B_{\frac{q-1}4}'\},$$
  where  $B_{i}=\{ H\overline{u^jl^i}\di j\in \ZZ_q\}$ and  $B_{i}'=\{ H\overline{tu^jl^i}\di j\in \ZZ_q\},$ where $1\le i\le \frac{q-1}4$.

\begin{lemma}
\label{sub}
Suppose $q\equiv 1({\rm mod}~ 4)$. Then for $G$ acting on $[G:H]$, \begin{enumerate}
\item [{\rm (1)}]\, there are $\frac{q-3}{2}$ suborbits of length $\frac{q+1}{2}$,  while  $\frac{q-1}{4}$ of them $\{H\overline{l^it}H\di 1\le i\le \frac{q-1}4\}$ are self-paired and  $\frac{q-5}{4}$ of them $\{H\overline{l^i}H\di 1\le i\le \frac{q-1}4\}$ are non-self-paired suborbits;
 \item [{\rm (2)}]\, there are $\frac{q-1}{4}$ suborbits of length $q+1$, with the form $H\overline{u^i}H$, where $i^2\in S^* \cap(4\th+N)$.  All of them are self-paired.
\end{enumerate}
\end{lemma}
\demo Since $q+1\equiv 3(\mod 4)$, for any $g\in G$, $H\cap H^g$ is either $\ZZ_2$ or 1, so every suborbit is of length either $\frac{q+1}2$ or $q+1$.

\vskip 3mm
 (1) $|\D|=\frac{q+1}2$

  \vskip 3mm Let $\Delta=HgH$ be  a  suborbit of length $\frac{q+1}{2}$. Then $H^{g}\cap H\cong Z_{2}$ and so  $\a^g$ is an involution of $H$,  where $\alpha=\overline{l^{\frac{q-1}4}}\in H$. Then
  $\a^g=\a^h$ for some $h\in H$, and so $gh^{-1}\in C_{G}(\alpha)=\lg \oll, \,\olt\rg $. Since $HgH=Hgh^{-1}H$, we may choose $h=1$ so that $g\in  C_{G}(\alpha)$.  Set
  $g=\overline{l^i}$ or $\overline{l^it}$ for some $i$.
 Moreover,  direct computations show that for any two distinct elements $g_1,\, g_2\in  C_{G}(\alpha)=\lg \oll,\, \olt\rg $, $Hg_1H=Hg_2H$ if and only if $g_1=g_2\a$. Therefore, we  have $\frac{q-1}2$ suborbits of length $\frac{q+1}2$.
 In particular, $HgH=Hg^{-1}H$ if and only if either $g^2=1$ or $g^{-1}=g\a$, where the second case gives $g\in H$. So we get $\frac{q-1}4$ self-paired suborbits $HgH$ where   $g$ is non-central involution in $C_{G}(\alpha)$ , noting $Hg\a H=HgH$. So the remaining  $\frac{q-5}{4}$ suborbits of length $\frac{q+1}{2}$ are non self-paired.

 \vskip 3mm  (2)  $|\D|=q+1$  \vskip 3mm

Let first consider the suborbits $D=H\overline{u^i}H$ where $i\in \ZZ_q^*$.   From the arguments in (1), we know that $|\D|=q+1$.
Since  $H\overline{u^{i}}H=H\a \overline{u^i} \a H=H\overline{u^{-i}}H$, $\D$ is self-paired. Set $g=\overline{u^i}$.

 Suppose that  $H^{g}\cap H=Z_{2}$, that is
 $$\overline{u^{-i}}\overline{t'(x_1,\,y_1)}\overline{u^{i}}\in H,$$
which implies $2x_1-iy_1=0$. Insetting it in $x_1^2-y_1^2\th=1$, we get $$i^2=4\th+4x_1^{-2}\in S^*\cap (4\th +S^*).$$ Therefore,  $\D$ is of length $q+1$ if and only if  $i^2 \in S^*\cap(4\th+N)$.
By Proposition \ref{pro:preseki}, $|S^*\cap(4\th+N)|=\frac{q-1}4$.
Check that $H\overline{u^i}H=H\overline{u^j}H$ if and only  if $i=\pm j$. Therefore, we get $\frac{q-1}4$ suborbits of length $q+1$.

\vskip 3mm
Since $1+\frac{q-3}2\frac{q+1}2+\frac{q-1}4 (q+1)=\frac{q(q-1)}2=|[G:H]|$, we already find all suborbits.\qqed

\vskip 3mm

In what follows we deal with all cases of  suborbits $\D$ in Lemma~\ref{sub}, separately.
\begin{lemma}
Suppose that $\D$ is self-paired and of length $\frac{q+1}2$. Then
$X(G,\,\D)$  is hamiltonian.
\end{lemma}
\demo
Let $X=X(G,\,\D)$, where $\D$ is self-paired and of length $\frac{q+1}2$. From the last lemma,  $\D=H\overline{l^kt}H$ for some $k$.
Note $\frac{q-1}4=r$ is a prime, the two smallest values for $q$  are 13 and 29.  One may find a $H$-cycle by Magma for $q=13$ and 29. So let $q\ne 13,\, 29$.
First we give a remark.

\vskip 3mm {\it Remark:} Suppose we  may get two facts: (i) for any $B'\in {\bf B'}$,  $d(H,\, B')=0, 2$ or 4;  (ii)  $d(H, \, \cup_{B'\in \B'}B')\ge 5$.
Then  $H$ is adjacent to at least two blocks $B'_i,\, B_j'$ in ${\bf B}'$ such that
 $d(H,\, B_i')=2$ or 4.  Let $Y$ be the block graph. Then $Y$ is a  bipartite graph   of order $2r$, where $r=\frac{q-1}4$ is a prime. Note that $H\in B_{r}$.
 Since $\lg \oll\rg /\lg \oll^{r}\rg $ acts regularly on both ${\bf B}$ and ${\bf B'}$, we may set $B_i'^{d}=B_j'$ for some $d\in \lg \oll\rg /\lg \oll^{r}\rg $.
  Then we get a $H$-cycle of $Y$:
   $$B_i',\, B_{r},\, B_i'^d,\, B_r^d,\, B_i'^{d^2},\, \cdots,\, B_r^{d^{r-1}},\, B_i'.$$
   Then by Proposition~\ref{cyclelift},  we may find a $H$-cycle for $X(G,\,\D)$.
\vskip 3mm
Come back to the proof. Check that the neighborhood of $H$ is:
 $$X_1(H)=\D=H\overline{l^kt}H=\{H\overline{l^kt}\overline{t(x_1,\,y_1)}\di x_{1}^2-y_{1}^2\th=1\}.$$
The vertex $H\overline{l^kt}\overline{t(x_1,\,y_1)}$  is contained in $B_i'$
 if and only if
 $$H\overline{l^kt}\overline{t(x_1,\,y_1)} = H\overline{tu^jl^i},\,{\rm for~some~}j,$$
 if and only if
$$\overline{l^kt}\overline{t(x_1,\,y_1)}(\overline{tu^jl^i})^{-1}\in H,$$
 if and only if  one of the following two systems of equations with unknowns $j$, $i$, $x_1$ and $y_1$  has solutions corresponding to $(\e,\, \eta)= (1,\, -1)$ and $(-1,\, 1)$:
 \begin{eqnarray} \left\{
\begin{array}{lll}
y_1j\th^{2k}&=&x_1(\th^{2k+2i}-\e),\\
y_1(\th^{2i+2}+\eta \th^{2k})&=&x_1\th j,\\
 x_1^2-y_1^2\th&=&1.
\end{array}
\right.  \end{eqnarray}
Every such  system has the same solutions with
\begin{eqnarray} \left\{
\begin{array}{llll}
 y_1^2&=&\frac{\th^{2k+1}}{\eta \th^{4k}+\th^2\e} \th^{2i}-\frac{\e \th}{\eta \th^{4k}+\th^2\e}, &\quad (i)\\
  y_1^2&=&\th^{-1}x_1^2-\th^{-1},  &\quad (ii)\\
 j&=&(\th^{2i}-\e \th^{-2k})\frac{x_1}{y_1}. &\quad (iii)
\end{array}
\right. \end{eqnarray}

From (iii), we know that given a solution for $x_1^2,\, y_1^2$ and $i$,  we have two values of $j$, that is $\pm j$.  Then the possible values for  $d(H,\,B_i')$ is 0, 2 or 4, noting we have two choices for $(\e,\, \eta)$, showing fact (i).

Set $b=-\th^{-1}$, $a_1=\frac{\th^{2k+1}}{\eta \th^{4k}+\th^2\e}$ and $a_2=-\frac{\e \th}{\eta \th^{4k}+\th^2\e}$. Then $a_1$, $a_2\neq 0$ and $a_2\ne b$.
From (i) and (ii), we get that either $$y_1^2\in S^*\cap (S^*+a_2)\cap (N+b)\,\, {\rm if}\,\, a_1\in S^*\quad {\rm or}\quad  y_1^2\in S^*\cap (N+a_2)\cap (N+b)\,\, {\rm if}\, \, a_1\in N.$$
By using  Lemma~\ref{feng}, we get that the number of solutions for  $y_1^2$ is at least  $\lfloor\frac{1}8(q-11-2\sqrt{q})\rfloor$, which implies that the  number of solutions for
$j,\, i,\, x_1,\,y_1$ is at least  $2\lfloor\frac{1}8(q-11-2\sqrt{q})\rfloor$, for given $(\e,\,\eta)$. In other words, $d(H,\, \cup_{B'\in \B'}B')$ is at least
$2\lfloor\frac{1}8(q-11-2\sqrt{q})\rfloor$. Moreover, $2\lfloor\frac{1}8(q-11-2\sqrt{q})\rfloor\ge 5$, showing fact (ii).
\qqed

\begin{lemma}
Suppose that $\D$ is non self-paired and of length $\frac{q+1}2$. Then
$X(G,\,\D\cup \D^*)$  is hamiltonian.

\end{lemma}
\demo
Let $X=X(G,\,\D\cup \D^*)$, where $\D$ is non self-paired and of length $\frac{q+1}2$. From Lemma~\ref{sub}, $\D=H\overline{l^k}H$ and $ \D^*=H\overline{l^{-k}}H$ for some integer $k$.
Note $\frac{q-1}4=r$ is a prime, the three smallest values for $q$  are 13, 29 and 53.  One may find a $H$-cycle by Magma for $q=13$, 29  and 53. So let $q\ne 13,\, 29,\, 53$.

From the remark in last lemma, it suffices to show   two facts: (i) for any $B'\in {\bf B'}$,  $d(H,\, B')=0,\, 2,\, 4,\, 6$ or 8;  (ii)  $d(H,\,  \cup_{B'\in \B'}B')\ge 9$.
\vskip 3mm

Check that the neighborhood of $H$ is:
$$X_1(H)=\D\cup \D^*=\{H\overline{l^k}\overline{t(x_1,\,y_1)},\, H\overline{l^{-k}}\overline{t(x_1,\,y_1)}\di x_{1}^2-y_{1}^2\th=1 \}.$$
The vertex $H\overline{l^k}\overline{t(x_1,\,y_1)}$ and $H\overline{l^{-k}}\overline{t(x_1,\,y_1)}$ are contained in $B_i'$ if and only if either

$$H\overline{l^k}\overline{t(x_1,\,y_1)}=H\overline{tu^jl^i},\ {\rm{or}}\, H\overline{l^{-k}}\overline{t(x_1,\,y_1)}=H\overline{tu^jl^i},\,{\rm{for\, some}}\, j$$
if and only if either
  $$\overline{l^k}\overline{t(x_1,\,y_1)}(\overline{tu^jl^i})^{-1}\in H,\ {\rm{or}}\, \overline{l^{-k}}\overline{t(x_1,\,y_1)}(\overline{tu^jl^i})^{-1}\in H$$
 if and only if one of the following four systems of equations with unknowns $j$, $i$, $x_1$ and $y_1$  has solutions corresponding to $(\e,\, \eta)= (1,\,-1),\,(1,\,1),\,(-1,\,-1)$ or $(-1,\,1)$:
  \begin{eqnarray} \left\{
\begin{array}{lll}
y_1(\th^{i+\epsilon k+1}-\eta \th^{-i-\e k})&=&x_1j\th^{\e k-i},\\
y_1j\th^{-\e k-i+1}&=&x_1(\th^{i-\epsilon k+1}-\eta \th^{-i+\e k}),\\
 x_1^2-y_1^2\th&=&1.
\end{array}
\right.  \end{eqnarray}
Every such system has the same solutions with
\begin{eqnarray} \left\{
\begin{array}{llll}
 y_1^2&=&\frac{\th^{2i}\th}{\eta \th\th^{2\e k}-\eta\th\th^{-2\e k}} -\frac{\eta \th^{2\e k}}{\eta \th\th^{2\e k}-\eta\th\th^{-2\e k}}, &\quad (i)\\
  y_1^2&=&\th^{-1}x_1^2-\th^{-1},  &\quad (ii)\\
 j&=&\frac{\th^{i+\e k+1}-\eta \th^{-i-\e k}}{\th^{\e k-i}}\frac{y_1}{x_1}. &\quad (iii)
\end{array}
\right. \end{eqnarray}

From (iii), we know that given a solution for $x_1^2,\, y_1^2$ and $i$,  we have two values of $j$, that is $\pm j$.  Then the possible values for  $d(H,\,B_i')$ is 0, 2, 4, 6 or 8, noting we have four choices for $(\e,\, \eta)$, showing fact (i).

Set $b=-\th^{-1}$, $a_1=\eta \th\th^{2\e k}-\eta\th\th^{-2\e k}$ and $a_2= -\frac{\eta \th^{2\e k}}{\eta \th\th^{2\e k}-\eta\th\th^{-2\e k}}$. Then $a_1$, $a_2\neq 0$ and $a_2\ne b$.
From (i) and (ii), we get that either $$y_1^2\in S^*\cap (N+a_2)\cap (N+b)\,\, {\rm if}\,\, a_1\in S^*\quad {\rm or}\quad  y_1^2\in S^*\cap (S^*+a_2)\cap (N+b)\,\, {\rm if}\, \, a_1\in N.$$
By using  Lemma~\ref{feng}, we get that the number of solutions for  $y_1^2$ is at least  $\lfloor\frac{1}8(q-11-2\sqrt{q})\rfloor$, which implies that the  number of solutions for
$j,\, i,\, x_1,\,y_1$ is at least  $2\lfloor\frac{1}8(q-11-2\sqrt{q})\rfloor$, for given $(\e,\,\eta)$. In other words, $d(H,\, \cup_{B'\in \B'}B')$ is at least
$2\lfloor\frac{1}8(q-11-2\sqrt{q})\rfloor$.
Moreover, $2\lfloor\frac{1}8(q-11-2\sqrt{q})\rfloor\ge 9$, showing fact (ii).
\qqed

\begin{lemma}
Suppose that $\D$ is self-paired and of length $q+1$. Then
$X(G,\,\D)$  is hamiltonian.

\end{lemma}
\demo
Let $X=X(G,\,\D)$, where $\D$ is self-paired and of length $q+1$. From lemma~\ref{sub}, $\D=H\overline{u^k}H$ for some integer $k$. Note $\frac{q-1}4=r$ is a prime.

 If we may  get two facts: (i) for any $B'\in {\bf B'}$,  $d(H,\, B')=0, 2$ or 4;  (ii)  for any $B'\in {\bf B'}$,  $d(H,\, B')=0, 2$ or 4,  then   every  vertex in block graph
 has the valency at least   $\frac{(q+1)-2}{4}=\frac{q-1}{4}=\frac 1 2 \frac{q-1}{2}$. So $Y$ contains a $H$-cycle. Since $d(B_i,\, B_j')$ is even,   this cycle can lift a $H$-cycle for $X(G,\,\D)$ by
 Proposition~\ref{cyclelift}.

\vskip 3mm

In fact, check that the neighborhood of $H$ is:
 $$X_1(H)=\D=H\overline{u^k}H=\{H\overline{u^k}\overline{t(x_1,\,y_1)}, H\overline{u^k}\overline{t'(x_1,\,y_1)}\di x_{1}^2-y_{1}^2\th=1 \}. $$By observing the neighbor, one can see these neighbors contained in $B_{\frac{q-1}{4}}$ are just $H\overline{u^k}$ and $H\overline{u^{-k}}$, which implies $d(B_{\frac{q-1}{4}})=2$.
The vertex $H\overline{u^k}\overline{t(x_1,\,y_1)}$ and $H\overline{u^k}\overline{t'(x_1,\,y_1)}$ are contained in $B_i$ if and only if either:
 $$H\overline{u^k}\overline{t(x_1,\,y_1)} = H\overline{u^{j}l^i},\,{\rm{or}}\,\ H\overline{u^k}\overline{t'(x_1,\,y_1)} = H\overline{u^{j}l^i}$$ if and only if either:
$$\overline{u^k}\overline{t(x_1,\,y_1)}(\overline{u^{j}l^i})^{-1}\in H,\,{\rm{or}}\,\ \overline{u^k}\overline{t'(x_1,\,y_1)}(\overline{u^{j}l^i})^{-1}\in H$$
 if and only if one of the following systems of equations with unknowns $j$, $i$, $x_1$ and $y_1$  has solutions corresponding to $(\epsilon,\,\eta,\,\gamma,\,\delta)=(-1,\,1,\,-1,\,1),\,(1,\,-1,\,1,\,-1),\,(-1,\,-1,\,-1,\,-1)$ or $(1,\,1,\,1,\,1)$ :
  \begin{eqnarray} \left\{
\begin{array}{lll}
\epsilon \th^{-i}y_1j&=&(x_1+ky_1)\th^{-i}-\eta x_1\th^{i},\\
\gamma(x_1+ky_1)\th^{-i}j&=&y_1\th^{-i}\th- \delta\th^{i}(y_1\th+kx_1),\\
 x_1^2-y_1^2\th&=&1.
\end{array}
\right.  \end{eqnarray}
This system has the same solutions with
$$j=\frac{(\th^{-i}-\eta\th^{i})x_1}{\e\th^{-i}y_1}+k\e^{-1}\,{\rm where}\, \d\e\th^{2i}=\gamma(k^{2}y_1^{2}+2kx_1y_1+1).$$
 Calculating the equation $\d\e\th^{2i}=\gamma(k^{2}y_1^{2}+2kx_1y_1+1)$ we could get
   $$(4k^2\th-k^4)u^2+(2k^2+2\d\e\gamma\th^{2i}k^2)u-(\d\e\th^{2i}-\gamma)^2=0,$$
where $u=y_1^2$. Since the product of the two solutions is $\frac{-(\d\e\th^{2i}-\gamma)^2}{4k^2\th-k^4}$, a non-square (as $4\th-k^2 \in N$), there exists at most one solution for $u=y_1^2$. It is easy to see that there are two solutions for $j$. Since there are just two different equations for $\d\e\th^{2i}=\gamma(k^{2}y_1^{2}+2kx_1y_1+1)$, there are at most $4$ solutions for $j$, that is $d(H,\,B_i)=0,\,2\,{\rm  or}\,4$, showing fact (i).

The vertex $H\overline{u^k}\overline{t(x_1,\,y_1)}$ and $H\overline{u^k}\overline{t'(x_1,\,y_1)}$ are contained in $B_i'$ if and only if either
$$H\overline{u^k}\overline{t(x_1,\,y_1)}= H\overline{tu^jl^i}\,{\rm{or}}\,\  H\overline{u^k}\overline{t'(x_1,\,y_1)}= H\overline{tu^jl^i}$$
 if and only if either
 $$\overline{u^k}\overline{t(x_1,\,y_1)}(\overline{tu^jl^i})^{-1}\in H \,{\rm{or}}\,\ \overline{u^k}\overline{t'(x_1,\,y_1)}(\overline{tu^jl^i})^{-1}\in H$$
 if and only if one of the following systems of equations with unknowns $j$, $i$, $x_1$ and $y_1$  has solutions corresponding to $(\epsilon,\,\eta,\,\gamma,\,\delta)=(1,\,-1,\,1,\,-1),\,(1,\,1,\,1,\,1),\,(-1,\,-1,\,-1,\,-1)$ or $(-1,\,1,\,-1,\,1)$:
\begin{eqnarray} \left\{
\begin{array}{lll}
-(x_1+ky_1)\th^{-i}j&=&\eta y_1\th^{-i}-\epsilon(y_1\th+kx_1)\th^{i-1}\\
-y_1\th^{-i}\th j&=&\delta(x_1+ky_1)\th^{-i}-\gamma x_1\th^{i-1}\th\\
x_1^2-y_1^2\th&=&1.
\end{array}
\right.  \end{eqnarray}
This system has the same solutions with
$$j=\frac{\d\th^{-i}-\gamma\th^{i}\th}{-\th^{-i}\th}\frac{x_1}{y_1}-\frac{k\d}{\th}\,{\rm where}\, \gamma\th^{2i}\th=\d(k^{2}y_1^{2}+2kx_1y_1+1).$$
Calculating the equation $\gamma\th^{2i}\th=\d(k^{2}y_1^{2}+2kx_1y_1+1)$ we could get
   $$(4k^2\th-k^4)u^2+(2k^2+2\d\gamma k^2\th^{2i}\th)u-(-\gamma\th^{2i}\th+\d)^2=0,$$
where $u=y_1^2$. Since the product of the two solutions is $\frac{-(-\gamma\th^{2i}\th+\d)^2}{4k^2\th-k^4}$, a non-square (as $4\th-k^2 \in N$), there exists at most one solution for $u=y_1^2$ and it is easy to see there are two solutions for $j$.
 Since there are just two different equations for $\gamma\th^{2i}=\d(k^{2}y_1^{2}+2kx_1y_1+1)$, there are at most $4$ solutions for $j$, that is $d(B_{\frac{q-1}{4}},\,B_i')=0,\,2\,{\rm  or}\,4$, showing fact (ii).
\qqed

\subsection{Groups in Table 2}
In this subsection, we  shall deal with the  groups in Table 2, separately.
\begin{lemma}
Let $G$ be one of  groups  in rows 1 and 2 of Table 2.
Then every  orbital graph of $G$ contains a Hamilton cycle.
\end{lemma}
\demo Let $T=\PSL(m,\,q)$ where $m=4$ or 5.  It suffices to consider the group  $T$. We shall deal with two cases: $m=4$ and $m=5$, separately.

  \vskip 3mm
 {\it Case 1:} $m=4$.
 \vskip 3mm
Let $\O$  be the set of 2-dim. subspaces of a space $V$ of dimension 4. Then  $n=\frac{(q^{4}-1)(q^{3}-1)}{(q-1)(q^{2}-1)}=(q^{2}+q +1)(q^{2}+1)$, where $s=q^2+q+1$ and $r=\frac{q^2+1}2$ are two   primes. Pick up  a subspace $W_0$ of dimension $d(W_0)=2$. Then  $T$ has two nontrivial suborbits  relative  to $W_0$:
 $$\D_1=\{W\in \O\di d(W\cap W_0)=1\} \quad {\rm and}\quad \D_2=\{W\in \O\di d(W\cap W_0)=0\},$$
 where
  $r_1:=|\D_1|=\frac{q^{4}-q}{q^{2}-q}=\frac{q^{3}-1}{q-1}$ and $r_2:=|\D_2|=n-1-r_1$.  Since $r_2\ge \frac n2$, the corresponding  orbital graph $\G(T,\, \D_2)$ has a H-cycle.

 Now  we are considering $X(T,\, \D_1)$.
 Take a projective  point $\lg \a\rg $ and   extend it into a base $\a,\, \a_1,\,\a_2,\,\a_3$ of $V$.
   Let $\si(\a)$ be the set of  all 2-dim. subspaces containing $\a$. Then $|\si(\a)|=q^2+q+1$.
 Since $\lg \a_1,\,\a_2,\,\a_3\rg $ contains exactly $q^2+q+1$ points and for any two distinct points $\b, \,\b'$ in $\lg \a_1,\,\a_2,\,\a_3\rg $, $\lg \a,\, \b\rg\ne \lg \a, \,\b'\rg $,
 one may see
 $$\si(\a)=\{ \lg\a,\, \b\rg \di \b\in \lg \a_1,\, \a_2, \a_3\rg \}.$$
  Let $\lg h\rg$ be the Singer subgroup of $\PSL(3,q)$ and $\beta\in\lg \a_1,\,\a_2,\,\a_3\rg$. Since $s=q^2+q+1$ is a prime,  $\lg \beta \rg,\, \lg \beta^{h}\rg ,\,\lg \beta^{h^{2}}\rg ,\,\cdots,\, \lg \beta^{h^{s-1}}\rg $ are
   all the projective points of $\lg \a_1,\,\a_2,\,\a_3\rg$. Denote $\beta^{h^i}=\beta_i$.
 Since the subgraph induced  by $\si (\a)$ is a complete graph, we may pick up a H-cycle of the subgraph, say
    $$\lg \alpha,\,\beta_{0}\rg,\, \lg \a,\, \beta_1\rg,\, \lg \alpha,\,\beta_{2}\rg ,\,\cdots, \,\lg \alpha,\,\beta_{s-2}\rg ,\, \lg\alpha,\,\beta_{s-1} \rg, \, \lg \alpha,\,\beta_{0}\rg ,$$
  where $\b_i\in \lg \a_1,\, \a_2, \,\a_3\rg $ and $s=q^2+q+1$.

   Set
  $$A=\{\lg \beta_i,\,\beta_{i+1} \rg,\,\lg\beta_{s-1},\,\beta_0 \rg|i=0,\,1,\,\cdots,\,s-2\}=\{ \lg \beta,\,\beta^{h}\rg ^{h^ {i}} |i=0,\,1,\,\cdots,\,s-1 \},  $$
  $$ X_i=\si(\b_i)\setminus ( \bigcup_{j=1}^{i-1}\si(\b_j)\bigcup A ), \, i\in \{ 1,\,2,\, \cdots, \, s-1\},$$
  $$ X_0=\si(\b_0)\setminus ( \bigcup_{j=1}^{s-1}\si(\b_j)\bigcup A ).$$
Since every 2-subspace $\lg\eta,\gamma\rg$  can be expressed as $\lg\eta,\, \gamma-\frac{b_0}{a_0}\eta\rg ,$ where  $\eta=a_0\a+a_1\a_1+a_2\a_2+a_3\a_3$ and $\gamma=b_0\a+b_1\a_1+b_2\a_2+b_3\a_3$,  every 2-subspace of $V$   is contained in $(\bigcup_{i=0}^{s-1} X_i) \bigcup A$. Moreover, from the definition, we know that $X_0,\, X_1,\, \cdots , \, X_{s-1}, \,A$ are mutually disjoint.

Now we are ready to  find a H-cycle for  $X(T, \,\D_1)$.   For $i=0,\,1,\,\cdots,\, r-2$ ,
   pick up a H-path $H_{i+1}$  in the subgraph induced  by $X_{i+1}\bigcup\{\lg \beta_i,\,\beta_{i+1} \rg \}$ with the starting vertex  $\lg \b_{i},\, \beta_{i+1}\rg$ and the ending vertex  $\lg \a, \,\beta_{i+1}\rg$.
   Pick up  $H_{0}$ in the subgraph induced  by $X_{0}\bigcup\{\lg \beta_{s-1},\,\beta_{0} \rg \}$ with the starting vertex $\lg \b_{s-1},\, \beta_{0}\rg$ and the ending vertex   $\lg \a,\, \beta_{0}\rg$.
   Then by replacing every arc ($\lg \alpha,\,\beta_i\rg , \lg\alpha,\,\beta_{i+1}\rg $) by  the path ($\lg \alpha,\,\beta_i\rg ,\, H_{i+1}$) and   the arc ($\lg \alpha,\,\beta_{s-1}\rg ,\, \lg\alpha,\,\beta_{0}\rg $) by the path ($\lg \alpha,\,\beta_{s-1}\rg ,\,    H_{0}$), we get a cycle:
   $$ \lg \alpha,\,\beta_{0}\rg,\,H_1,\,H_2,\, \cdots ,\,H_{s-1},\,H_0,  $$
   which is clearly   a H-cycle of $X(T,\, \D_1)$, as shown in Figure 1.
   
   \begin{figure}[!htp]
\begin{center}
\includegraphics[width=0.9\textwidth]{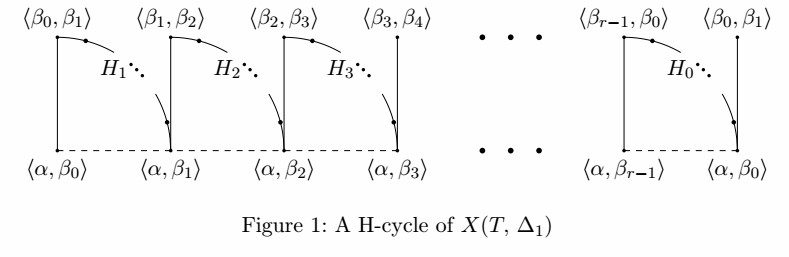}
\end{center}
\end{figure}
\include{pictue}

 {\it Case 2:} $m=5$.
 \vskip 3mm
Let $\O$ be the set of 2-dim. subspaces of $V$. Then
$$n=|\O|=\frac{(q^{5}-1)(q^{4}-1)}{(q-1)(q^{2}-1)}=(q^{4}+\cdots +1)(q^{2}+1)=2rs.$$
Then $s=q^{4}+\cdots +1$  is a prime  and $r=\frac{q^2+1}2$ are two prime.  Let $S= \lg h\rg$ be a Singer subgroup of $\PSL(5,q)$, where   $|S|=s$.
Take a projective point $\a$. Then $\alpha,\,\alpha^{h},\,\cdots,\,\alpha^{h^{s-1}}$ are all the projective points.  Set $W_{i}=\lg\alpha,\,\alpha^{h^ i}\rg$ where $i=1,\,2,\,\cdots,\,s-1$.
Then  $G$ has two nontrivial suborbits   relative  to $W_1$:
$$\D_1=\{W\in \O|d(W \cap W_{1})=1\}\quad  {\rm and}\quad \D_2=\{W\in \O|d(W \cap W_{1})=0\},$$
 where
  $$\begin{array}{ll}&r_1:=|\D_1|=(\frac{q^{4}}{q-1}-1)(q+1)=q(q+1)(q^2+q+1),\\
  &r_2:=|\D_2|=\frac{(q^5-q^2)(q^5-q^3)}{(q^2-1)(q^2-q)}=q^{4}(q^2+q+1).\end{array}$$
    Since $r_2\ge \frac n2$, the corresponding  orbital graph $X(T,\, \D_2)$ has a H-cycle.

 Now  we are considering $X(T,\,\D_1)$.
 Let $S_i$ be the path
$${W_{i},\,W_{i}^{h^{i}},\,W_{i}^{h^{2i}},\,W_{i}^{h^{3i}},\,\cdots,\, W_{i}^{h^{(s-1)i}}}.$$
Since $\lg h^i \rg $ acts nontrivially on $W_i$ and it is of order a prime $s$, $\lg h^i\rg $ moves $W_i$.
Since every 2-subspace must  be contained in some  clique and either $|S_i\cap S_j|=0$ or $S_i= S_j$ for any two distinct cliques $S_i$ and $S_j$, we could pick up $q^2+1$ distinct cliques which cover all 2-dim. subspaces, denoted by $W_{\mu_1},\,W_{\mu_2},\,\cdots,\, W_{\mu_{q^2+1}}$.
Then we can get  a H-cycle of $X(T,\,\D_1):$
$$W_{\mu_1},\,W_{\mu_1}^{h^{\mu_1}},\, W_{\mu_1}^{h^{2\mu_1}},\,\cdots, W_{\mu_1}^{h^{(s-1)\mu_1}},\,W_{\mu_2},\,W_{\mu_2}^{h^{\mu_2}},\,W_{\mu_2}^{h^{2\mu_2}},\,W_{\mu_2}^{h^{3\mu_2}},\,\cdots, W_{\mu_{q^2+1}}^{h^{(s-1)\mu_{q^2+1}}},\, W_{\mu_1}.$$
\qqed

\begin{lemma}
Every orbital graph of $G={\rm P\O}^{-}(2m,\,q)$ in row 3 of Table 2 is hamiltonian.
\end{lemma}
\demo Let $G={\rm P\O} ^{-}(2m,q)$ act on $n$ totally singular $1$-spaces, where $n=\frac{(q^{m}+1)(q^{m-1}-1)}{q-1}=2rs$ and $m=2^{2^l}$.  Then $m-1$ is a prime. Since $m-1=(2^{2^{l-1}}-1)(2^{2^{l-1}}+1)$,  we get $2^{2^{l-1}}-1=1$, which implies  $l=1$ and then $m=4$. Now $r=:\frac{q^{3}-1}{q-1}$ is  a prime.
 Let $\Omega$ be the set of all t.s.$1$-spaces.
Recall that $\SO^{-}(8,\,q)\leq \GL(8,\,q)$ and $|\GL(8,\,q)|=q^{28}\Pi_{i=1}^8(q^i-1)$.
To describe $\SO ^-(8,\,q)$, take  a symmetric bilinear form,  given by the following matrix:
$$J=\left(
              \begin{array}{ccc}
                0 & E_{3} & 0 \\
                 E_{3} & 0 & 0 \\
                0 & 0 & J_{2} \\
              \end{array}
            \right),\quad
J_{2}=\left(
        \begin{array}{cc}
          1 & 0 \\
          0 & -t \\
        \end{array}
      \right), \quad t\in N.$$
Let $\langle A\rangle$ be a Singer subgroup of $\GL(3,\,q)$, $C=A^{q-1}$ and  $D=(C^{-1})'$,  where $C'$ denotes the transpose of $C$.
 Set $B=C\bigoplus (C')^{-1}\bigoplus E_2$, the block diagonal matrix.   Then we have $BJB'=J$, which means $B\in \SO ^{-}(8,\,q)$. Since $\overline{B}$ is of prime order, $\overline{B}\in (\PSO^{-}(8,\,q))'={\rm P\O}^{-}(8,\,q).$
  Set $S=\lg \overline{B}\rg $ and $\alpha=(1,\,0,\,\cdots,\,0)$. Then there are two nontrivial  suborbits for the action of $G_{\lg\alpha\rg}$ relative to $\lg \a\rg$,  see \cite{LWWX94}:
$$\D_1=\{ \lg \b\rg \in \O\setminus\{\lg a\rg\} \di (\a,\, \b)=0\}  \quad {\rm and}\quad  \D_2=\{ \lg \b\rg \in \O\setminus \{\lg a\rg\} \di (\a,\, \b)\ne 0\}, $$
where  $|\Delta_{1}|=q^5+q^4+q^2+q$ and  $|\Delta_{2}|=q^{6}$. Since $|\Delta_{2}|\geq \frac{1}{2}n$, we only need to consider $X(G, \,\D_1)$.

    Noting that  $S$ acts semiregularly on $\Omega$, we consider  the block graph $\ox$ induced by   $S$-orbits,  where $V(\ox)=q^4+1$.
For any $\gamma =(\g_1,\, \g_2,\, \g_3)\in \O$, where $\g_1=(c_1,\, c_2,\, c_3)$, $\g_2=(c_4,\, c_5,\, c_6)$ and $\g_3=(c_7,\, c_8)$,  we have  $\g \overline{B}^i J \a'=0$ if and only if   $\g_2D^i \a'=0$, that is
$\g_2D^i=(0,\, c_5',\, c_6')$ for some $c_5',\, c_6'$. Since $\lg C\rg$  (and so $\lg D\rg $) is regular on nonzero $1$-spaces,  we know that
  $\alpha$ has $q+1$ (resp. $q^2+q$) neighbors in the block  $\gamma ^S$ if $\gamma \not\in \a^S$ (resp. $\g\in \a^S$).  From
  $((q^5+q^4+q^2+q)-(q^2+q))/(q+1)=q^4$ we know that $\ox$ is a complete graph. By Propsosition~\ref{cyclelift}, $X(G,\, \D_1)$ is  hamiltonian.

\begin{lemma}
Every orbital graph of $G={\rm P\Omega}^{+}(2m,\,q)$ in row 4 of Table 2 is hamiltonian.
\end{lemma}
\demo
Let $G={\rm P\Omega}^{+}(2m,\,q)$ act on $n$ totally singular $1$-spaces, where $n=\frac{(q^{m}-1)(q^{m-1}+1)}{q-1}=2rs$,   $m=2^{2^l}+1$,  and  $s=\frac{q^{m}-1}{q-1}$ and $r=\frac{q^{m-1}+1}2$ are primes. Let $\Omega$ be the set of all totally singular $1$-spaces.
Recall that ${\rm SO}^{+}(2m,\,q)\leq \GL(2m,\,q)$.
To describe $\SO ^+(2m,\,q)$, take  a symmetric bilinear form,  given by the following matrix:
$$J=\left(
              \begin{array}{cc}
                0 & E_{m}  \\
                 E_{m} & 0  \\

              \end{array}
            \right).
$$

Let $\langle A\rangle$ be a Singer subgroup of $\GL(m,\,q)$, $C=A^{q-1}$ and  $D=(C^{-1})'$,  where $C'$ denotes the transpose of $C$.
Set $B=C\bigoplus (C')^{-1}$. Then we have $BJB'=J$, which means $B\in \SO ^{+}(2m,\,q)$.
Since $B$ is of prime order, $\overline{B}\in (\PSO^{+}(m,q))'={\rm P\O}^{+}(m,\,q).$
Set $S=\lg \overline{B}\rg $  and
 $\alpha=(1,\,0,\,\cdots,\,0)$. Then there are two nontrivial  suborbits for the action of $G_{\langle\alpha\rangle}$ relative to $\langle\alpha\rangle$, see By \cite{LWWX94}:
$$\D_1=\{ \lg \b\rg \in \O\setminus\{\lg a\rg\} \di (\a,\, \b)=0\}  \quad {\rm and}\quad  \D_2=\{ \lg \b\rg \in \O\setminus \{\lg a\rg\} \di (\a,\, \b)\ne 0\}, $$
where
$|\Delta_{1}|=\frac{(q^{m-1}+q)(q^{m-1}-1)}{q-1}$ and $|\Delta_{2}|=q^{2m-2}.$ Since $|\Delta_{2}|\geq \frac{1}{2}n$, we only need to consider $X(G, \,\D_1)$.

Noting that  $S$ acts semiregularly on $\Omega$, we consider  the block graph $\ox$ induced by   $S$-orbits,  where $V(\ox)=q^{m-1}+1$.
For any $\gamma =(\g_1,\, \g_2)\in \O$,  we have  $\g \ols^i J \a'=0$ if and only if   $\g_2D^i \a'=0$, which implies that the first coordinate of $\g_{2}D^i$ is 0.
 Since $\lg C\rg$  (and so $\lg D\rg $) is regular on nonzero $1$-spaces,  we know that
  $\alpha$ has $\frac{q^{m-1}-1}{q-1}$ (resp. $\frac{q^{m}-1}{q-1}-1$) neighbors in the block  $\gamma ^S$ if $\gamma \not\in \a^S$ (resp. $\g\in \a^S$).  From
  $(\frac{(q^{m-1}+q)(q^{m-1}-1)}{q-1}-(\frac{q^m-1}{q-1}-1))/(\frac{q^{m-1}-1}{q-1})=q^{m-1}$ we know that $\ox$ is a complete graph. By Propsosition~\ref{cyclelift}, $X(G,\, \D_1)$ is  hamiltonian.\qqed

\begin{lemma}
Vertex-transitive graphs arising from the action of $A_{c}$ on $2$-subsets given in row 6 of Table 2 are hamiltonian.
\end{lemma}
\demo
Let $\Omega=\{\alpha_{1},\,\alpha_{2},\,\cdots,\,\alpha_{c}\}$, where $c\ge 5$.  Then  we only have the following two orbital graphs:

 \vskip 3mm (1) Two subsets are adjacent if and only if they intersect  at a single point.  Then we may get a $H$-cycle as the following way:

  first pick up  a cycle of  $c$ vertices, say $\{\alpha_{1},\,\alpha_{2}\},\, \{\alpha_{2},\,\alpha_{3}\}, \, \{\alpha_{3},\,\alpha_{4}\},\, \cdots,\, \{ \a_{c-1},\,\alpha_{c}\},\, \{\alpha_{c},\,\alpha_{1}\},$ \quad $\{\alpha_{1},\,\alpha_{2}\};$
 then

replace the edge  $\{\alpha_{1},\, \alpha_{2}\},\, \{\alpha_{2},\,\alpha_{3}\}$ by    any $H$-path of  all 2-subsets  containing $\a_{2}$, with the starting vertex   $\{\alpha_1,\, \alpha_2\}$  and the ending vertex  $\{\alpha_{2},\,\alpha_{3}\}$; then

  replace the  edge $\{\alpha_{2}, \,\alpha_{3}\},\, \{\alpha_{3},\,\alpha_{4}\}$ by  any $H$-path of   all 2-subsets  containing $\a_{3}$,  with the starting vertex   $\{\alpha_2,\, \alpha_3\}$  and the ending vertex  $\{\alpha_{3},\,\alpha_{4}\}$; then for $5\le i\le c$,

replace the  edge $\{\alpha_{i-2},\, \alpha_{i-1}\}, \,\{\alpha_{i-1},\,\alpha_{i}\}$ by  any $H$-path of   all 2-subsets  containing $\a_{i-1}$ but removing  $\{ \{\a_2,\, \a_{i-1}\}, \,  \{\a_3,\, \a_{i-1}\}, \, \cdots,  \, \{\a_{i-3}, \,\a_{i-1}\} \},$  with the starting vertex   $\{\alpha_{i-2},\,\alpha_{i-1}\}$  and the ending vertex  $\{\alpha_{i-1},\,\alpha_{i}\}$.

\vskip 3mm
 (2)  Two subsets are adjacent if and only if they have no intersecting point. In this case,    if $c\ge 7$, then the degree of the graph is more than  $\frac n2$ and so it  is hamiltonian. For $c\le 5$, we do it just by Magma.
\qqed

\begin{lemma} Let $G$ be one of the  groups  listed in row 5, 7-10 of Table 2. Then every orbital graph of $G$
is hamiltonian.
\end{lemma}
\demo Using Magma, we get  compute the suborbits for these groups and show that every corresponding orbital graph is hamiltonian.

(1) The action of $\PSL(3,5)$ on the flags  has  three nontrivial  suborbits, with the respective  length   10, 50 and 125;

(2) The action of $M_{11}$ on the cosets of a subgroup isomorphic to $S_{5}$ has three nontrivial suborbits, with the respective  length  15, 20 and 30;

(3) The action of $M_{12}$ on the cosets of a subgroup isomorphic to $M_{10}:2$ has two nontrivial suborbits, with the respective  length $20$ and $45$;

(4) The action of $M_{23}$ on the cosets of a subgroup isomorphic to $A_{8}$ has three nontrivial  suborbits, with the respective  length  $15$, $210$ and $280$;

(5) The action of $J_{1}$ on the cosets of a subgroup isomorphic to $\PSL(2,11)$ has four nontrivial  suborbits, with the respective  length   $11$, $12$, $110$ and $132$.
\qqed

\vskip 3mm
\f {\bf Acknowledgments:}
This work is partially   supported by
the National Natural Science Foundation of China (12071312 and 11971248). All authors declare that  this paper has  no conflict of interest.

{\footnotesize

\end{document}